%babushkas.tex: In How Many Ways Can You Reassemble Several Russian Dolls?
%%a Plain TeX file by and Doron Zeilberger (8 pages)

%begin macros
\def\D{{\cal D}}

\def\P{{\cal P}}
\def\Q{{\cal Q}}
\baselineskip=14pt
\parskip=10pt

\font\eightrm=cmr8 
\font\eighttt=cmtt8
\magnification=\magstephalf
\def\P{{\cal P}}
\def\Q{{\cal Q}}
\def\1{{\overline{1}}}
\def\2{{\overline{2}}}
\parindent=0pt
\overfullrule=0in

\def\frac#1#2{{#1 \over #2}}
%\headline={\rm  \ifodd\pageno  \RightHead  \else  \LeftHead  \fi}
%\def\RightHead{\centerline{
%Title
%}}
%\def\LeftHead{ \centerline{Doron Zeilberger}}
%end macros
\bf
\centerline
{
In How Many Ways Can You Reassemble Several Russian Dolls?
}
\rm
\bigskip
\centerline{ {\it Doron ZEILBERGER}\footnote{$^1$}
{\eightrm  \raggedright
Department of Mathematics, Rutgers University (New Brunswick),
Hill Center-Busch Campus, 110 Frelinghuysen Rd., Piscataway,
NJ 08854-8019, USA.
%\break
{\eighttt zeilberg at math dot rutgers dot edu} ,
\hfill \break
{\eighttt http://www.math.rutgers.edu/\~{}zeilberg} .
Sept. 16, 2009.
Accompanied by Maple package {\eighttt BABUSHKAS}
downloadable from {\eighttt http://www.math.rutgers.edu/\~{}zeilberg/tokhniot/BABUSHKAS} .
Exclusively published in the Personal Journal of Ekhad and Zeilberger
{\eighttt http://www.math.rutgers.edu/\~{}zeilberg/pj.html} and {\eighttt http://arxiv.org} .
Supported in part by the USA National Science Foundation.
}
}

\quad
{\it Dedicated to my favorite identical twins: Thotsaporn and  Thotsaphon THANATIPANONDA}

{\bf Preface}

In the Fall of 1981, I had the pleasure and honor of living next-door to
the eminent combinatorialist Joel Spencer, who was then a visiting professor at the Weizmann Institute
of Science in Israel. Joel
was always glad to talk to me, provided that I spoke in Hebrew. Joel is also a great punner,
and was very proud when he made his first pun in Hebrew, during a volleyball game with other
math faculty and students, when he said ``kadur sheli'' that could mean ``my ball'' 
but also ``Lee [Segal]'s ball''.
One day I got as a present a set of ``Russian Dolls'' (a.k.a. as {\it Matryoshkas} or {\it Babushkas}),
which is a nested set of dolls.
Trying to test Joel's knowledge of enumeration (after all he is an expert in
``Hungarian'', rather than enumerative, combinatorics), I asked him:

{\it In how many ways can one reassemble an $n$-nested Russian Doll?},

and he immediately replied: {\it This is a Stirling question, it rings a Bell}.

{\bf Several Russian Dolls}

Almost thirty years later, my brilliant student, Thotsaporn ``Aek'' Thanatipanonda
asked me what happens if you have {\it several}, say $r$, identical Russian Dolls.
Thanks to {\it Sloane}, {\it MathSciNet}, and {\it Google scholar}, we
quickly found out that the case $r=2$  goes back to Comtet[C]
and is featured in Sloane's {\bf A020554}.
See also [Ba] and [R], and for an insightful {\it species} treatment see [L] and [P].
The case $r=3$ was nicely handled by Ed Bender[Be], while the
general case was given its {\it coup de gr\^ace} by John Devitt and David Jackson[DJ],
who used a very ingenious generatingfunctionology approach.

As Herb Wilf[W1] famously said, an ``answer'' to an enumeration question is an {\it efficient algorithm} to generate 
many terms in the enumerating sequence. Explicit formulas are just {\it one} such way, and often
not the most efficient one! In this article, I will describe a ``calculus'' approach, using
{\it differential operators}, that has the following advantages.

{\bf 1.} It is somewhat faster than using the Devitt-Jackson[DJ] complicated exponential generating function.

{\bf 2.} It is so flexible that it can handle {\it non-identical} Russian Dolls.
Given
$a_1$ single-births, $a_2$ pairs of (identical) twins, $a_3$ sets of (identical) triplets, $a_4$ sets of (identical) quadruplets, 
$\quad \dots \quad$, $a_k$ sets of (identical) $k$-tuplets, $\dots$, (and you can't tell identical twins etc. apart form each other),
in how many ways can you partition them into (not necessarily distinct) {\it sets}?
This is useful if the school principal has to
assign them into (non-ordered) classes such that no identical-looking children would be in the
same class (or else their teacher won't be able to tell them apart, and the children can play tricks on her).

{\bf 3.} It can be used in conjunction with Wilf's [Wi2] celebrated methodology for random selection
of combinatorial objects to design quick algorithms for selecting, {\it uniformly at random},
a mutli-set set-partition of any given multiset.

{\bf 4.} It beautifully illustrates MacMahon's
 lovely  method of ``differential operators'' that transcribes combinatorial operations into differential operators.
MacMahon was very  fond of it,
and he nicely described it in the entry {\it Combinatorial Analysis} of the eleventh edition of {\it Encyclopeadia Britannica}
[M].

{\bf 5.} It beautifully illustrates my favorite methodology of {\it rigorous experimental mathematics}.
You teach the computer how to do the combinatorics, it derives, {\it all by itself}, the (symbolic) differential operator
(that for $r>3$ would be too complicated to derive by hand, let-alone apply, even for MacMahon), and then the computer
goes on and uses it to crank-out as many terms as desired in the enumerating sequence.

{\bf 6.} It beautifully illustrates the notion of {\it catalytic variables}. These are variables
corresponding to quantities that we may not care about, but are nevertheless needed in order to
facilitate the enumeration. At the end of the day, we set them all equal to $1$.

{\bf 7.} Last but not least, it enables me to 
contribute six new sequences to {\it Sloane}, the (beginnings of) the enumerating sequence for the number of
ways of reassembling $r$ identical $n$-nested Russian Dolls for $3 \leq r \leq 8$.
So far, only $r=1$ (the Bell numbers, {\bf A000110}) and $r=2$ (Comtet's sequence {\bf A020554}) are present there.

{\bf The Evolution Differential operator}

{\it A Baby Example}: $r=1$

For the sake of pedagogy, let's first treat the classical case of one Russian Doll.

Suppose that you have a {\it set partition} of $\{1,2, \dots, n-1\}$ with $k$ sets. 
How can we accommodate a new-comer $n$? Either she is shy (or anti-social), and
decides to form her own new set $\{n\}$, or she is outgoing and would like to join
an existing set, for which she has $k$ choices. If we call the ``number of sets''
the ``state'', then in the former case the new comer, $n$, caused the set-partition
to move to state $k+1$, while in each of the $k$ latter cases, it stayed in
state $k$. If we give state $k$ the {\it weight} $z^k$, then each and every
set-partition of state $k$ (with weight $z^k$) gives rise to the ``evolution''
$$
z^k \rightarrow z^{k+1}+ k z^k \quad  .
$$
In other words
$$
z^k \rightarrow (z+ z{{d} \over {dz}}) z^k  \quad .
$$
This is true for each and every monomial $z^k$, and by linearity for each polynomial.
So if $P_n(z)$ is the sum of the weights of all set-partitions of $\{1, \dots, n\}$, then
we have the {\it differential-recurrence} equation:
$$
P_{n}(z)= \D_1 P_{n-1}(z) \quad ,
$$
where $\D_1$ is the {\it differential operator}
$$
D_1 f(z) := (z+ z{{d} \over {dz}})f(z) \quad .
$$
The initial condition is $P_{0}(z)=1$. If we are only interested in
the {\it total number} of set partitions of $\{1,2, \dots , n\}$, then at the
end of the day we plug-in $z=1$, getting
$$
B_n =P_n(1) \quad .
$$
This gives a quick way to crank-out a table of the first one thousand (or whatever) Bell numbers,
that is memory efficient. Once we are at day $n$, and know $P_{n}(z)$, we can let the computer
forget about $P_{n-1}(z)$ (and even about $B_{n-1}$, once it is printed out).
Of course, this is not more efficient than using
$$
B_n=\sum_{i=0}^{n-1} {{n-1} \choose {i}} B_i \quad ,
$$
or
$$
\sum_{n=0}^{\infty} {{B_n} \over {n!}} z^n= e^{e^z-1} \quad ,
$$
but (probably) not less efficient, and besides we only did it as a warm-up.

{\it A Toddler Example}: $r=2$ 

Suppose that we have $n-1$ pairs of identical twins, already arranged into sets.
The only restriction is that no two identical twins can be in the same set,
but sets can be repeated. So we have a {\it multiset of sets}, whose union, as
a multiset, is $1^2 2^2 \dots (n-1)^2$.

Of course, no set can be repeated more than twice (why?).
Let there be $a$ sets, $A_1, \dots, A_a$ that show up {\it once}, and
let there be $b$ sets $B_1, \dots , B_b$ that show up twice.
We say that the {\it state} of this arrangement is $(a,b)$, and its {\it weight} is
$z_1^az_2^b$.

The school principal has to place the newly-arrived pair of twins $n$ and $n$,
but he can't put them in the same set. There are {\bf seven} cases.

{\bf Case 1}: Create two new singleton sets $\{n \}^2$. There is only {\it one} way of
doing it, and the new state is $(a,b+1)$, with weight $z_1^a z_2^{b+1}$.
The corresponding operation on monomials is
$$
f(z_1,z_2) \rightarrow z_2 f(z_1,z_2) \quad .
$$

{\bf Case 2}: Create one new singleton set $\{n \}$,
and place the other $n$ into one of the existing $A_i's$.
There are {\it a} ways of
doing it, and the new state is $(a+1,b)$, with weight $z_1^{a+1} z_2^b$.
The corresponding operation on monomials is
$$
f(z_1,z_2) \rightarrow z_1 \left ( z_1 {{d} \over {d z_1}} \right ) f(z_1,z_2) \quad .
$$

{\bf Case 3}: Create one new singleton set $\{n \}$,
and place the other $n$ into one of the existing $B_i's$.
There are {\it b} ways of
doing it, and the new state is $(a+3,b-1)$, 
since one of the $B$'s became two $A$s, with the help of $n$.
The weight of the new state is $z_1^{a+3} z_2^{b-1}$.
The corresponding operation on monomials is
$$
f(z_1,z_2) \rightarrow z_1^3 z_2^{-1} \left ( z_2 { {d} \over {d z_2}} \right ) f(z_1,z_2) \quad .
$$

{\bf Case 4}: Place the two twins $n$ and $n$ into
two (different, of course) $A_i$'s.
There are ${{a} \choose {2}}$  ways of
doing it, and the new state remains $(a,b)$, with weight $z_1^a z_2^b$.
The corresponding operation on monomials is
$$
f(z_1,z_2) \rightarrow {{1} \over {2}} \left ( z_1^2 {{d^2} \over {d ^2z_1}} \right ) f(z_1,z_2) \quad .
$$

{\bf Case 5}: Place one of the new twins ($n$ or $n$) into
one of the $A_i$'s, and the other one into one of the $B_i$'s.
There are $ab$  ways of
doing it, and the new state is $(a+2,b-1)$,
since one of the doubletons $B$'s,
let's call it $B_i$, was lost, and it became the two distinct sets
$B_i$ and $B_i \cup \{ n\}$.
The new weight is $z_1^{a+2} z_2^{b-1}$.
The corresponding operation on monomials is
$$
f(z_1,z_2) \rightarrow z_1^2 z_2^{-1}   \left ( z_1 {{d} \over {z_1}} \right )
\left ( z_2 {{d} \over {dz_2}} \right ) f(z_1,z_2) \quad .
$$

{\bf Case 6}: Place the twins ($n$ and $n$) into
two different $B_i$'s.
There are ${{b} \choose {2}}$  ways of
doing it, and the new state is $(a+4,b-2)$,
since two of the doubletons $B$'s,
let's call them $B_i$ and $B_j$, were lost, and they became the four distinct new sets
$B_i$, $B_i \cup \{ n\}$, $B_j$ and $B_j \cup \{ n\}$.
The new weight is $z_1^{a+4} z_2^{b-2}$.
The corresponding operation on monomials is
$$
f(z_1,z_2) \rightarrow z_1^4 z_2^{-2}     {{1} \over {2}} \left ( z_2^2 {{d^2} \over {dz_2^2}} \right ) f(z_1,z_2) \quad .
$$

{\bf Case 7}: Place both of the new twins ($n$ and $n$) into
each of the two copies of the same  $B_i$.
There are $b$  ways of
doing it, and the new state is the same, $(a,b)$,
since single sets stay single and doubletons stay doubletons.
The new weight is $z_1^a z_2^b$.
The corresponding operation on monomials is
$$
f(z_1,z_2) \rightarrow  \left ( z_2 {{d} \over {dz_2}} \right ) f(z_1,z_2) \quad .
$$

Combing, we have just proved:

{\bf Fact}: Let $P^{(2)}_n(z_1, z_2)$ be the sum of the weights of all multiset set-partitions of
the multiset $\{ 1^2 \dots n^2 \}$, with the weight being
$z_1$ to the power the number of sets that show up once times $z_2$ to the power the number of sets that show
up twice. Let $\D_2$ be the partial-differential operator
(where $D_1:={{d} \over {dz_1}}$, $D_2:={{d} \over {dz_2}}$),
$$
\D_2 :=
z_{{2}}D_{{2}}+1/2\,{z_{{1}}}^{4}{D_{{2}}}^{2}+{z_{{1}}}^{3}D_{{1}}D_{{2}}+1/2\,{z_{{1}}}^{2}{D_{{1}}}^{2}+{z_{{1}}}^{3}D_{{2}}+{z_{{1}}}^{2}D_{{1}}+z_{{
2}} \quad ,
$$
then
$$
P^{(2)}_n(z_1,z_2)=\D_2 P^{(2)}_{n-1}(z_1,z_2) \quad .
$$
The Comtet numbers are $P^{(2)}_n(1,1)$.

{\bf The general case}

Suppose that we already have a multiset set-partition of $\{1^{r} \dots (n-1)^{r}\}$,
with $a_1$ sets that show-up once, $a_2$ sets that show-up twice, $\dots$, $a_r$ sets
that show-up $r$ times. The state of this particular multiset set-partion is
$(a_1,a_2, \dots, a_r)$, and its weight is $z_1^{a_1} \dots z_r^{a_r}$.

We have to place the $r$ identical new comers $n^r$.

We must make the following decisions

{\bf 1.} How many of them would start their own singleton sets, say, $c_0$

{\bf 2.} For the remaining $r-c_0$ new members $n$, for $i=1,2, \dots , r$,
how many of them, let's call it $c_i$, would be placed in sets that show up
$i$ times. 

After these decisions we have a vector of non-negative integers 
$[c_0,c_1, \dots, c_r]$, such that $c_0+c_1+ \dots + c_r=r$.

Once we decided that $c_i$ of the new $n$'s would go to sets that show up $i$ times,
we have to decide, amongst those $c_i$ siblings which sets should be asked to
invite them. These sets can be all different, but they could all be different copies of
the same set (if $c_i \leq i$). This naturally leads to an {\it integer partition}
$\lambda_i=1^{m_1} 2^{m_2} 3^{m_3}\dots i^{m_i}$ (written in multiplicity notation).

We have to place $m_1$ of these siblings such that each of them goes to different sets.
These make $m_1$ of the formerly $a_i$-repeated sets become $(a_i-1)$-repeated sets
and creates $m_1$ new singleton sets. So $a_1 \rightarrow a_1+m_1$, and $a_{i-1} \rightarrow a_{i-1}+m_1$
and $a_i \rightarrow a_i -m_1$. In terms of differential operators it is
$$
f(z_1, \dots , z_r) \rightarrow (1/m_1!) (z_1 z_{i-1} D_{i})^{m_1} f(z_1, \dots , z_r)  \quad .
$$
Similarly for $2^{m_2}$ we have
$$
f(z_1, \dots , z_r) \rightarrow (1/m_2!) (z_2 z_{i-2} D_{i})^{m_2} f(z_1, \dots , z_r)\quad ,
$$
and so on.

So every possible scenario of placing the new identical $r$ siblings of the $n$ family corresponds
to an r+1 tuple:

$$
T=[c_0, \lambda_1, \dots, \lambda_r] \quad ,
$$
where $c_0$ is an integer, $\lambda_1, \lambda_2, \dots, \lambda_r$ are integer partitions such
that the largest part of $\lambda_i$ is $\leq i$, and
$$
c_0+ |\lambda_1|+ |\lambda_2|+ \dots + |\lambda_r|=r \quad .
$$
For each such scenario corresponds the ``monomial'' operator
$$
\P[T]:=z_{c_0} \prod_{i=1}^r \Q[\lambda_i] \quad ,
$$
where, writing $\lambda_i=1^{m_1} 2^{m_2} 3^{m_3}\dots i^{m_i}$ ($m_1,m_2$ are now {\it local variables}, i.e.
they  are different, of course, for each $\lambda_i$),
$$
\Q[\lambda_i]=\prod_{j=1}^{i} {{1} \over {m_j!}} (z_jz_{i-j}D_{i})^{m_j} \quad .
$$
Here $D_i:={{d} \over {dz_i}}$ and $z_0:=1$.

Finally, we can write down the {\it evolution operator} $\D_r$:
$$
\D_r :=\sum_{T \,\, scenario} \P[T] \quad,
$$
and we have the

{\bf Theorem}: Let $P^{(r)}_n(z_1, z_2, \dots, z_r)$ be the sum of the weights of all multiset set-partions of
the multiset $\{ 1^r \dots n^r \}$, with the weight being
$z_1$ to the power the number of sets that show up once times $z_2$ to the power the number of sets that show
up twice times $\dots$ times $z_r$ to the power the number of sets that show up $r$ times, then
$$
P^{(r)}_n(z_1,z_2, \dots , z_r)=\D_r P^{(r)}_{n-1}(z_1,z_2, \dots , z_r) \quad .
$$
The {\it number} of such multi-set set partitions is of course
$P^{(r)}_n(1,1, \dots , 1)$.

{\bf Non-Identical Russian Dolls}

The operators $\D_r$ can be combined to yield the

{\bf Main Theorem}: The number of ways of partitioning into sets
a multiset consisting of $m_1$ elements that appear once,
$m_2$ elements that appear twice, $\dots$, $m_r$ elements that appear $r$ times,
or equivalently, the multiset
$$
1 \dots m_1 (m_1+1)^2 \dots (m_1+m_2)^2 \dots (m_1+ \dots + m_{r-1}+1)^{r} \dots  (m_1+ \dots + m_{r-1}+m_r)^{r}
$$
is computed as follows.
First compute the polynomial in $z_1, \dots, z_r$:
$$
P(z_1,z_2, \dots , z_r)=\left ( \prod_{i=1}^{r} \D_{i}^{m_i} \right ) (1) \quad ,
$$
and then plug-in $z_1=1, \dots, z_r=1$.

{\bf Random Generation} 

Going back to the Stirling-Bell case, $r=1$, the 
``differential recurrence'' $P_n(z)=(z+ z{{d} \over {dz}}) P_{n-1}(z)$ is
equivalent to the famous recurrence
$$
S(n,k)=S(n-1,k-1)+kS(n-1,k) \quad .
$$
This can be used, according to Wilf[W2], to generate {\it uniformly at random},
a set-partion with $k$ sets as follows.

First pre-compute a table of $S(n,k)$ using the recurrence. Now roll a loaded coin with
probability of Heads being $S(n-1,k-1)/S(n,k)$ and probability of Tails being
$kS(n-1,k)/S(n,k)$. If it lands Heads, 
recursively generate a random set-partitions of $\{1, 2, \dots, n-1\}$ with $k-1$ sets, and adjoin the singleton
$\{n\}$ to it, otherwise  generate
recursively a random set-partition of $\{1, 2, \dots, n-1\}$ with $k$ sets, 
and then roll a fair $k$-sided die,
and accordingly decide which of the $k$ members of the set-partition should invite $n$ to join it.

If we want a (uniformly) random set partition, then decide on the number of sets $k$, by
rolling a loaded $n$-faced die with probabilities
of it landing $k$ equalling $S(n,k)/B_n$, and then proceed as before.

The differential-recurrence of the Theorem yields to a partial recurrence for
the quantity, let's call it $S^{(r)}(n; a_1, \dots, a_r)$ for the number of
multiset set-partitions of $1^r \dots n^r$ with $a_1$ sets that show up once, $\dots$, 
$a_r$ sets that show up $r$ times. Using this 
the computer (all by itself!) can use the Wilf Methodology to create
a random-generation algorithm. The programming details are a bit daunting, so we leave it as a challenge
to the reader.

{\bf The Maple package BABUSHKAS}

Everything here (except for the random-generation, for which we only have the simple $r=1$ case)
is implemented in the Maple package {\tt BABUSHKAS} available, via a link, from the
webpage of this article:
{\eighttt http://www.math.rutgers.edu/\~{}zeilberg/mamarim/mamarimhtml/babushkas.html} \hfill\break
or directly from: {\eighttt http://www.math.rutgers.edu/\~{}zeilberg/tokhniot/BABUSHKAS} .
That webpage also contains sample input and output, including the sequences for $1 \leq r \leq 8$.

The main procedure is {\tt SeqBrn(r,n)} that uses the present approach to generate the first $n$ terms
of the enumerating sequence for the number of ways of reassembling $r$ identical Russian Dolls.
{\tt SeqBrnDJ(r,n)} does the same thing using the Devitt-Jackson approach. We are glad to report
that they agree! As yet another check, we have the program {\tt SSP}, that actually constructs
the set of {\it all} multi-set set partitions of any given multiset, and emables checking, for
small values, with the naive count. Procedure {\tt B(L)} handles the case
of non-identical Russian Dolls, or equivalently, an arbitrary multiset, using the
Main Theorem.

{\tt SeqCrn} and {\tt SeqCrnDJ} handle set-partitions of the multiset $1^r \dots n^r$, in other words,
each set can only show up once. This is simply $P^{(r)}_n(1, 0 , \dots , 0)$, in the
above notation.

Full details are available on-line by typing {\tt ezra();}.

{\bf The sequences}

Even though much more data is available in the above-mentioned webpage, and these
sequences will soon be submitted to {\it Sloane}, let us cite the first ten terms
for $r=1,2,3,4$.

$r=1: 1, 2, 5, 15, 52, 203, 877, 4140, 21147, 115975$ (the Bell Numbers).

$r=2: 1, 3, 16, 139, 1750, 29388, 624889, 16255738, 504717929, 18353177160$ (the Comtet numbers)

$r=3: 1, 4, 39, 862, 35775, 2406208, 238773109, 32867762616, 6009498859909, 1412846181645855$

$r=4: 1, 5, 81, 4079, 507549, 127126912, 55643064708, 38715666455777, 40095856807088486$, \hfill\break
$58901884724160709571$.

{\bf References}

[Ba] G. Baroti, {\it Calcul des nombres de bicouvrements et de birev\^etements d'un ensemble fini,
employant la m\'ethode fonctionnelle de Rota}, 
in: Combinatorial theory and its Applications I (ed. P. Erd{\H o}s et al.), North Holland, Amsterdam,
1970), 93-103.

[Be] E.A. Bender, {\it Partitions of multisets}, Discrete Mathematics {\bf 9}(1974), 301-312.

[C] L. Comtet, {\it Birecouvrements et birevetements d'un ensemble fini}, Studia Sci. Math. Hungar. {\bf 3}(1968), 137-152.

[DJ] J.S. Devitt and D.M. Jackson, {\it The enumeration of covers of a finite set}, J. London Math. Soc.(2) {\bf 25} (1982), 1-6.

[L] G. Labelle, {\it Counting enriched multigraphs}, Discrete Math. {\bf 217}(2000), 237-248.

[M] P.A. MacMahon, {\it Combinatory Analysis}, in:
{\it Encyclopeadia Britantica, 11th edition} (available on-line)
vol. 6, 752-758 (1910).

[P] G. Paquin, {\it D\'enombrement de multigraphes enrichis}, M\'emoire, Math. Dept., Univ. Qu\'ebec \`a Montr\'eal, 2004. 

[R] J. Reilly, {\it ``Bicoverings of a Finite Set'' by Generating Function Methods}, J. Comb. Theory, Series A {\bf 28} (1980), 219-225.

[W1] H.S. Wilf, {\it What is an answer?}, Amer. Math. Monthly {\bf 89} (1982), 289-292. 

[W2] H. S. Wilf,
{\it A unified setting for sequencing, ranking, and selection algorithms for combinatorial objects},
Advances in Mathematics {\bf 24}(1977), 281-291.
\end